\documentclass[12pt,twoside,english,a4paper]{article}
\usepackage{etex}
\usepackage{caption}
\usepackage[latin1]{inputenc}
\usepackage[intlimits] {amsmath}   % muss for defjs# stehen
\usepackage{graphicx}
\usepackage{psfrag}
\usepackage{ifthen}
\usepackage{fancyhdr}
\usepackage{rotating}
\usepackage{multirow}
\usepackage{booktabs}              % Dicke Tabellenlinien
\usepackage{amssymb}
\usepackage{amsbsy}
\usepackage{bm}                    % fette Schrift in Formeln
\usepackage{bbm}
\usepackage{babel}
\usepackage{theorem}
\usepackage{upgreek}  % gerade (roman) griechische Buchstaben z.B. \upalpha
\usepackage{pifont}   % zusaetzliche Schriften, z.B. eingekreiste Zahlen
\usepackage[active]{srcltx}   % aktive suche im xdv

\usepackage{cite}
\usepackage{url}
\makeatletter
\AtBeginDocument{%
  \def\doi#1{\url{https://doi.org/#1}}}
\makeatother
\usepackage{suetterl}  % Suetterlin Schrift
\newfont{\suetdbl}{suet14 scaled 2000}  % Suetterlin skaliert um Faktor 2,000

\usepackage{yfonts}  % altdeutsche Schrifte und Initialen
\newfont{\gothdbl}{ygoth scaled 2000} 
\newfont{\frakdbl}{yfrak scaled 2000}
\newfont{\swabdbl}{yswab scaled 2000}
\usepackage[normalem]{ulem}  % zusaetzliche Unterstreichungsformen wie z.B.:
\usepackage[textsize=footnotesize,textwidth=1.7cm]{todonotes}
\usepackage[]{defjs3}      % eigene Definitionen: Option BoldVarsRoman

\usepackage[authoryear]{natbib}
\bibpunct{[}{]}{,}{a}{}{;}
\fancypagestyle{plain}{%
  \fancyhead{}%
  \fancyfoot[c]{\sffamily\thepage}%
}
\makeatletter                           % Leere Fuellseiten
\def\cleardoublepage{\clearpage\if@twoside \ifodd\c@page\else
  \hbox{}
  \vspace*{\fill}
  \thispagestyle{empty}
  \newpage
  \if@twocolumn\hbox{}\newpage\fi\fi\fi}
\makeatother
\begin{document}
\unitlength1.0cm
\frenchspacing

%=== cover page
%\input{ftemplate_js_2}

%=== sections ===
%--- section 1

\captionsetup[figure]{labelfont={bf},labelformat={default},labelsep=period,name={Fig.}}
\captionsetup[table]{labelfont={bf},labelformat={default},labelsep=period,name={Tab.}}

\thispagestyle{empty}
\ce{\bf \large
A surrogate model for data-driven magnetic}

\ce{\bf \large stray field calculations}

\vspace{4mm}
\ce{Rainer Niekamp$^1$, Johanna Niemann$^1$, Maximilian Reichel$^1$, Hongbin Zhang$^2$, J\"org Schr\"oder$^1$}

\vspace{4mm}
\ce{$^1$Institute of Mechanics, University of Duisburg-Essen,}
\ce{Universit\"atsstr. 15, 45141 Essen, Germany}
\ce{\small e-mail: j.schroeder@uni-due.de,
    phone: +49 201 183 2708,
    fax: +49 201 183 2680}
\vspace{2mm}
\ce{$^2$Theory of Magnetic Materials, Technical University Darmstadt,}
\ce{Otto-Berndt-Straße 3, 64287 Darmstadt, Germany}
\ce{\small e-mail: hongbin.zhang@tu-darmstadt.de,
    phone: +49 6151 16-23135}

\vspace{4mm}
\begin{center}
{\bf \large Abstract}
\bigskip

{\footnotesize
\begin{minipage}{14.5cm}
\noindent
In this contribution we propose a data-driven surrogate model for the prediction of magnetic stray fields in two-dimensional random micro-heterogeneous materials.
Since data driven models require thousands of training data sets, FEM simulations appear to be too time consuming. Hence, a stochastic model based on Brownian motion, which utilizes an efficient evaluation of stochastic transition matrices, is applied for the training data generation.
For the encoding of the microstructure and the optimization of the surrogate model, two architectures are compared, i.e. the so-called UResNet model and the Fourier Convolutional neural network (FCNN). Here we analyze two FCNNs, one based on the discrete cosine transformation and one based on the complex-valued discrete Fourier transformation. Finally, we compare the magnetic stray fields for independent microstructures (not used in the training set) with results from the FE$^2$ method, a numerical homogenization scheme, to demonstrate the efficiency of the proposed surrogate model.

\end{minipage}
}
\end{center}


\section{Introduction}
Deep learning algorithms have successfully been employed in numerous applications that perform different tasks such as object recognition, approximation, optimization, classification, regression, and forecasting.
The fundamental framework of deep learning, the neural network, was already outlined in 1943, see \cite{McCPit:1943:tbo}. However, it has gained increasing importance and scientific interest due to Big Data approaches and growing computational capabilities.
From various successful applications involving Deep Learning, it is evident that this method is applicable within a large number of different disciplines.
DNNs have also been succesfully used for the fast determination of quasi-optimal coupling constraints for the FETI-DP domain decomposition solver \cite{HeiKlaLan:2021:cml},as fast surrogate models in computational fluid dynamics (CFD), see \cite{GuoLiuOer:2016:cnn,EicHeiKla:2020:scn} to name but a few.\\
Within this work, the simulation of demagnetization and stray magnetic fields was analyzed under the assumption of a time-independent linear problem with homogeneous magnetic material in free space and a finite domain. 
Solving magnetic stray fields can be challenging and time-consuming.
Deep learning models have the main advantage of providing a computationally cheaper option to traditional finite element frameworks for the computation.
The computation is particularly difficult with respect to the extension to infinite domains, which are typically associated with an open boundary value problem. 
Specific methods exist to approximate these open boundary value problems, such as FEM-based truncation methods, the numerical boundary element method, or the scaled boundary finite element method (SBFEM), see \cite{SchReiBir:2022:aen,BirReiSch:2022:msw,BucRucRai:2003:cbd,SonWol:1997:thb}.
However, these numerical methods are highly computationally intensive.
Therefore, the development of machine learning approaches for simulations of magnetic stray fields, especially in the infinite domain, may be of particular interest in the future and could provide an efficient and cost-effective alternative for the computation of magnetic stray fields.
In recent years, deep neural networks (DNNs) have also emerged in the scientific field of magnetism and magnetic materials.
In \cite{AdlAbd:2014:unn}, a DNN is applied for modeling complex magnetic materials and magnetic field calculations, where in \cite{NguBolBer:2022:emo} it is used for the prediction of magnetic field distributions generated by permanent magnets. 
Furthermore, Deep Convolutional Neural Networks (DCNNs), specifically a convolutional autoencoder, have been successfully used in micromagnetism to approximate the time evolution of the magnetization in a thin film, see \cite{KovFisOez:2019:lmd}. 
Employment of Physics Informed Neural Networks (PINNs) has been demonstrated for problems of magnetostatics, micromagnetics, and hysteresis computation in \cite{KovExlKor:2022:mam}, by using a deep neural network to approximate the magnetic vector potential.
The implemented U-shaped residual neural network (UResNet) and the Fourier Convolutional Neural Network (FCNN) each represent a different extension of Deep Convolutional Neural Networks.
Initially introduced by \cite{RonFisBro:2015:ucn} for applications in medical image segmentation and processing, the UResNet has since been a widely used deep learning technique in medical research, see \cite{RonFisBro:2015:ucn,HuaLinTon:2020:uaf,GuaKhaSik:2019:fdu}. 
Moreover, this type of segmentation network can be applied in other fields, such as e.g. fluid dynamics for the prediction of air quality, see \cite{JurReiBen:2022:dlm}, or for the prediction of steady-state laminar flow, see \cite{RibRehAhm:2020:ess,TakNas:2020:add}. A significant improvement for the training of very deep neural networks was the use of Residual Convolutional Neural Networks (ResNets), which have been previously applied, e.g., in the inverse design of magnetic fields, see \cite{PolBjoJor:2021:ido}. 
Furthermore, Fourier Convolutional Neural Networks (FCNNs) are another efficient enhanced method introduced by \cite{MatHenLeC:2014:fto} and described in \cite{PraWilCoe:2017:Ffc} with the approach of computing convolutions in the frequence domain to reduce computational cost and to accelerate training. In this manner, the Fourier Neural Operator (FNO), introduced by \cite{LiKovAzi:2020:fno}, yields an accurate solution of parametric partial differential equations such as Burgers' equation, Darcy flow, and Navier-Stokes equation. In \cite{WenLiAzi:2021:Uae}, a U-FNO was presented as an extension of the traditional FNO, combining the spectral domain convolution and U-net approaches for modeling highly complex multiphase flow problems.
In \cite{BelBilAli:2022:pin}, about 100 different network architectures of neural networks were compared for solving parametric magnetostatic problems. The predictive accuracy for each of these studies was compared to the results from finite element analysis. The overall result showed that incorporating a Fourier encoding in the frequency space can improve the co-generalization rate of the predictions. \\

In this contribution, two data-driven surrogate models built on Convolutional Neural Networks (CNNs) are presented to solve magnetostatic problems, particularly for calculating demagnetization and stray fields for two-dimensional geometries simulated with Brownian motion. 
First, the constitutive laws of electromagnetism will be presented in section \ref{sec:Magnetism} and gradually reduced for application in Magnetostatics.  
The three fundamental magnetic quantities, magnetic induction, magnetization, and magnetic field strength, are briefly discussed to set a contextual framework for the magnetic stray fields. 
Section \ref{sec:data} then outlines two methods to generate data, forming the foundation for the training optimization of the machine learning models.
 The computation of the magnetic fields can be performed by the computationally-intensive FEM or the less expensive stochastic Brownian motion.
 The stochastic Brownian motion, which is algebraically implemented using stochastic transition matrices, is a more efficient approach that allows the generation of a large number of data sets within a relatively short time.
 Prior to generating the new dataset, a cross-validation of the calculations from the FEM and Brownian motion simulations is performed.
In section \ref{UResNet} the first surrogate model, a modified U-shaped residual network (UResNet) is proposed, combining features and spatial information using a hierarchical structure.
The second surrogate model is presented in section \ref{FCNN} and is based on a Fourier Convolutional Neural Network (FCNN), where convolutions are computed in the Fourier domain by applying the complex-valued discrete Fourier transform (DFT) instead of in the spatial domain.
Within the second alternative approach, convolution in the spatial domain equals an element-wise product in the Fourier domain, which can significantly reduce the computational cost of the training.
The comparison of these two machine learning models with respect to their training characteristics and approximation performance using error measurements, a numerical example, and the correlations is eventually drawn in section \ref{sec:Example}.

\section{Theoretical Background of Magnetostatics}\label{sec:Magnetism}
Magnetostatics provides the physical foundation for describing magnetism in solids independent of time. It entails the study of magnetic fields and forces. Magnetic fields can be generated by electric currents or by the dispersion of magnetic material. The focus of this chapter is the outline of the theoretical basis of Maxwell's equations as well as the physical background and laws of the main magnetic fields.
\subsection{The Magnetic Induction $\mathbf{B}$}
In magnetism, Gauss's law for magnetic fields defines the magnetic induction, also known as magnetic flux density, to have a divergence equal to zero with
\begin{equation}
	\div \bB = 0\quad\textrm{and}\quad\bB = \hat{\bB}(\bH) \, .
	\label{eq:gauss_law}
\end{equation}
This implies that the magnetic field is solenoidal as the lines of force form continuous loops and no magnetic monopoles exist, see \cite{Poo:2018:glf}.
Due to the absence of magnetic charges, there is no net flux across a closed surface, see \cite{Jac:1999:ce}.
$\bB$ can be described as a function of the magnetic field strength denoted by $\bH$. 
As mentioned in \cite{Coe:2010:mag}, the following three main sources contribute to magnetic induction:
\begin{enumerate}
	\item electric currents in conductors
	\item moving charges (which generate electrical current)
	\item magnetic moments (equivalent to current loops).
\end{enumerate} 
Generally, the behavior of magnetic materials is assumed to be non-linear. However, the free space surrounding the material is assumed to have the linear relation
\begin{equation}
	\bB = \mu_0 \bH	\, .
	\label{Eq:free_space}
\end{equation}
When a magnet is placed in the free space characterized by (\ref{Eq:free_space}), a magnetization $\bM$ occurs inside the magnet and stray or demagnetizing fields develop. 
Within magnetizable media, magnetic induction is thus defined as the sum of magnetic field strength $\bH$ and the magnetization $\bM$ multiplied by the magnetic permeability of free space $\mu_0$
\begin{equation}
	\bB = \mu_0 (\bH + \bM) \, .
	\label{Eq:mu(H+M)}	
\end{equation}
This equation relates the three fundamental magnetic quantities. 
Fig. \ref{fig:mag_fields} illustrates $\bB$, $\bM$, and $\bH$ for a uniformly magnetized material without any externally applied magnetic field.  
\begin{Figure}[hbtp!]
	\centering
	\quad \includegraphics[width=14cm]{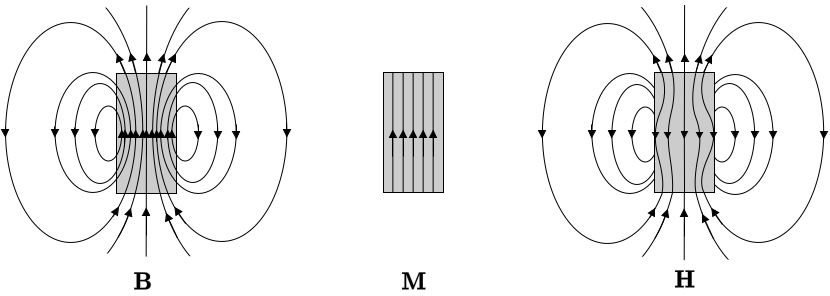}
	\caption{Illustration of the magnetic induction $\bB$, the magnetization $\bM$ and the magnetic field strength $\bH$, compare \cite{Leh:2018:dgd}}
	\label{fig:mag_fields}
\end{Figure}

\subsection{The Magnetization $\mathbf{M}$}
Magnetization refers to the impact magnetic material has on the magnetic induction given an existing external magnetic field.
Inside the magnet, a magnetization $\bM$ arises, resulting in the generation of a magnetic dipole, see \cite{Jil:1990:itm}. 
Due to the magnetic dipole, sources (sinks) appear at the surface of the magnet, which is contradictory to the Gauss Law for Magnetic Fields.
The sources (sinks) of magnetization $\bM$ must thus be considered sources (sinks) of magnetic field strength $\bH$. 
Mathematically, this is shown by
\begin{equation}
	\div \bB = \div \mu_0(\bH +\bM) = 0 
	\quad \Rightarrow \quad
	\div \bM = -\div \bH.
\end{equation}
Magnetization in so-called linear media, such as diamagnetic and paramagnetic materials, can be determined according to \cite{Coe:2010:mag,Jil:1990:itm} by
\begin{equation}
	\bM = \chi \bH
	\label{Eq:Chi_H}
\end{equation}
since these materials have constant magnetic susceptibilities $\chi$.
If no magnetization exists, this would lead back to the previous relation between $\bB$ and $\bH$ in free space. By substituting (\ref{Eq:Chi_H}) into (\ref{Eq:mu(H+M)}), we obtain 
\begin{equation}
	\bB = \mu_0 \bH + \mu_0 \bM = \mu_0 \bH +\mu_0 \chi \bH = \mu_0 (1+\chi) \bH
\end{equation}
with the relative permeability
\begin{equation}
	\mu_r = (1+\chi)
\end{equation}
and the permeability of the material
\begin{equation}
	\Bmu = \mu_0 \Bmu_r \, ,
\end{equation}
resulting in
\begin{equation}
	\bB = \Bmu \bH \, . 
\end{equation}

\subsection{The Magnetic Field Strength $\mathbf{H}$}
Each magnet generates a magnetic field in the surrounding space and within its volume. 
Due to this reason, the magnetic field strength is separated into an external and an internal field, see \cite{Coe:2010:mag}. The magnetic field $\bH$ can be split into
\begin{equation}
	\bH = \overline\bH+\widetilde\bH
\end{equation}
with
\begin{enumerate}
	\item $\overline\bH$ created by conduction currents or other external sources and
	\item $\widetilde\bH$ created by the magnetic solid. 
\end{enumerate}
The first contribution is the spatially constant external field, also called the externally applied field, to which the magnetic particle has no contribution.
The stray field occurring outside the magnet or the demagnetizing field that arises inside the magnet makes up the second component. 
The magnetic field strength $\widetilde{\bH}$ and the magnetization $\bM$ are oriented in opposing directions within the magnetic sample. Thus, $\widetilde{\bH}$ is called the demagnetizing field.
\section{Generation of Data with Stochastic Model}\label{sec:data}
Let $\B$ be the body of interest parameterized in $x$, and let the surface of $\B$ be denoted by $\partial\B$. 
We divide the latter into a part $\partial\B_B$, describing the magnetic flux density across the surface, and a remaining part $\partial\B_{\varphi}$, indicating the magnetic potential. 
The edge decomposition satisfies the relations $\partial\B=\partial\B_B\cup\partial\B_{\varphi}$ and $\partial\B_\varphi\cap \partial\B_B=\emptyset$. 
By the Gauss's law for magnetic fields, the differential equation underlying the magnetostatic boundary value problem is given as  
\begin{equation}
	\div \bB = 0 \quad \textrm{with} \quad \bB = \Bmu \bH,
	\label{Eq:BVP}
\end{equation}
with $\bB$ being the magnetic induction, $\bH$ the auxiliary field, and $\Bmu$ the magnetic permeability. As the curl of the gradient of any scalar-valued field is zero, this implies that we can introduce a scalar potential with
\begin{equation}
	\bH(\bx):= -\nabla\varphi,
	\label{Eq:magnetic_potential}
\end{equation}
with $\nabla(\boldsymbol{\cdot}) = \partial_{\bx}(\boldsymbol{\cdot})$. The boundary conditions can be defined in terms of the magnetic potential $\varphi_0$ and the magnetic flux $\zeta_0$ as 
\begin{equation}
	\varphi=\varphi_0\quad\textrm{on}\;\partial\B_\varphi
\quad\textrm{and}\quad	
	\zeta_0=\bB\cdot\bn \quad\textrm{on}\;\partial\B_B.	
	\label{Eq:BC}
\end{equation}
\subsection{Brownian Motion}
The modeling of Brownian motion can be implemented by using the Monte Carlo simulation introduced by \cite{MetRosRos:1953:eos} or by using stochastic matrices.
 The stochastic matrix approach is a method from the probability theory that characterizes the transition probabilities of Markov chains and thus enables the prediction of possible future distributions. 
Here, the simulation with stochastic transition matrices was selected since this approach represents a faster algebraic algorithm.
\subsubsection*{Stochastic Matrix Approach}
The simulation of Brownian motion with stochastic transition matrices is described using a representative heterogeneous microstructure with the geometric dimensions $L \cdot H$.
The microstructure within the given geometric region consists of two phases: phase A represents the vacuum space, and phase B is the magnetic material.
Let this region comprise a total number of $S$ positions defined by
\begin{equation}
	S=\{(i,j) \in [1,L] \times [1,H]\}
\end{equation}
where $i$ and $j$ denote the indices of the individual positions. 
For two adjacent positions within this region, $J=(i,j)$ denotes the current position and $I=(i,j)$ the neighboring position. 
At the initial state, the particles are randomly located at the left boundary $(j=1)$ with the same probability of distribution for each position. Thus, the initial state yields
\begin{equation}
     p_{\,0\,I}=\left\{\begin{array}{cl}
        m^{-1}, & \text{if } I=(i,1) \in S\\.
        0, & \text{else } \\
        \end{array}\right\}
       \label{eq:Initial_state}
\end{equation}
where $m$ is the total number of admissible positions on the left boundary.
Hereby noting that the sum of all initial, transition and remaining probabilities is 1 under all conditions. 
In other words, all particles are initially placed on the boundary, and after placement, no further option exists than for the particle to move to a new position or to remain in the current position. 
When positioned at the left boundary, the motion of each particle is subject to the so-called transition rules given below, which describe the permissible movements of each particle for each time step with the help of transition probabilities $t_{ij}$. 
\begin{enumerate}
\item The transition probability of the movement from a position $J$ in phase A to a new position $I$ in phase B is given by $t_{ij}=\frac{\mu_r^{A} \mu_r^{B}}{2(\mu_r^{A}+\mu_r^{B})}$, corresponding to $\frac{1}{4}$ of the harmonic mean. The probability of remaining in the same position $J$ is $1 -\sum{t_{ij}}$.
\item When moving within the same phase, one of the four direct neighbors $I$ (left, right, bottom, top) can be visited with probability $t_{ij}=\frac{\mu_r^{A}}{4}$ for phase A and $t_{ij}=\frac{\mu_r^{B}}{4}$ for phase B. 
\item If the particle wants to cross a vertical boundary ($j=1$ or $j=L$), the particle is moved to a randomly chosen admissible point on the left boundary
\item If the particle wants to cross a horizontal boundary edge ($i=1$ or $i=H$), the particle movement is rejected, and the particle stays at its current position. 
\end{enumerate}
These four transition rules used throughout our simulation are depicted in Fig. \ref{fig:Random_Walk}. 
\begin{Figure}[hbt]
	\centering
	\includegraphics[width=12cm]{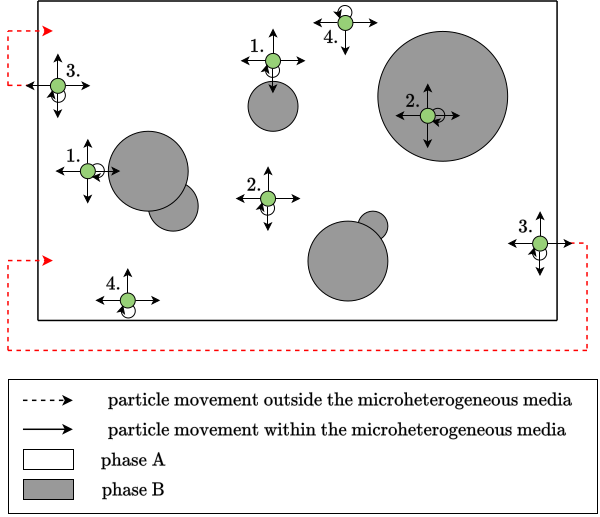}
	\caption{A visualization of the transition rules $1.$ to $4.$ and movement of the particles}
	\label{fig:Random_Walk}
\end{Figure}

These transition probabilities find application in the computation of the stochastic transition matrices, denoted as $\mathsf{M}$. 
The stochastic transition matrix $\mathsf{M}_{I,J} \in [0,1]$ comprises the probabilities of a particle moving from the current position $J$ to one of the neighboring positions $I$.
In general, the probability distribution of the next time step $n+1$ can be calculated from the multiplication of the transition matrix $\mathsf{M} \in \mathbb{R}^{[L\cdot H \times L\cdot H]}$ by the probability distribution from the current time step $n$ through
\begin{equation}
	\bp_{n+1} = \mathsf{M} \, \bp_{n}
\end{equation}
using the initial distribution from (\ref{eq:Initial_state}). We can perform this forward iteration as it allows us to compute a steady state since the transition matrix $\mathsf{M}$ has a left eigenvector with eigenvalue 1. 
Thus, if the eigenvectors of $\mathsf{M}$ are known, so is the stationary distribution defined as
\begin{equation}
\mathsf{M} \, \overset{*}\bp =\overset{*}\bp
\end{equation}
Let us now consider the example of a micro-heterogeneous material with the material parameters
$\mu_0 = 4 \pi \cdot 10^{-7} \,\,\text{H/m} \,\,\text{in A}$ and
	$\mu_1 = 2 \pi \cdot 10^{-7} \,\,\text{H/m} \,\,\text{in B}$
into consideration. The respective relative magnetic permeabilities are $\mu_r^{A} = 1$ in A and $\mu_r^{B} = 0.5$ in B.
For $\mu_r > 1$, the two relative permeabilities are rescaled by division through the larger value of the two material parameters $\mu_r^{\text{max}}$.
The transition probability that a particle located at position $J$ in phase A moves to position $I$ in phase B is given by
\begin{equation*}
	t_{ij}= \dfrac{\mu_r^{A}\cdot \mu_r^{B}}{2 \left(\mu_r^{A}+\mu_r^{B}\right)} = \frac{1}{6}
\end{equation*}
and the transition probabilities within a phase are calculated to be 
\begin{equation*}
	t_{ij} = \dfrac{\mu_r^{A}}{4} = \frac{1}{4} \quad
	\text{in A} \quad \text{and} \quad t_{ij} = \frac{\mu_r^{B}}{4} = \frac{1}{8} \quad
	\text{in B.}
\end{equation*}
Under the assumption of these material parameters and transition probabilities, the transition matrix $\mathsf{M}$ is then calculated as shown in Fig. \ref{fig:Transition_Matrix}. 
\begin{Figure}[hbt!]
	\centering
	\raisebox{-0.5\height}{\includegraphics[width=4.3cm]{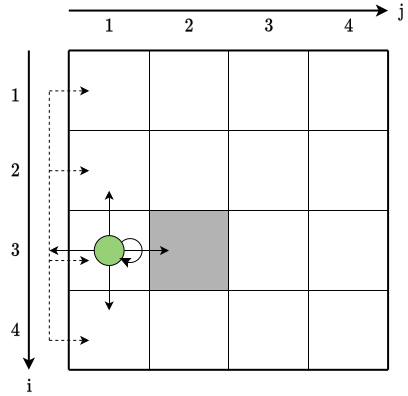}}
    \hspace{0.2cm}
	\raisebox{-0.5\height}{\includegraphics[width=11.2cm]{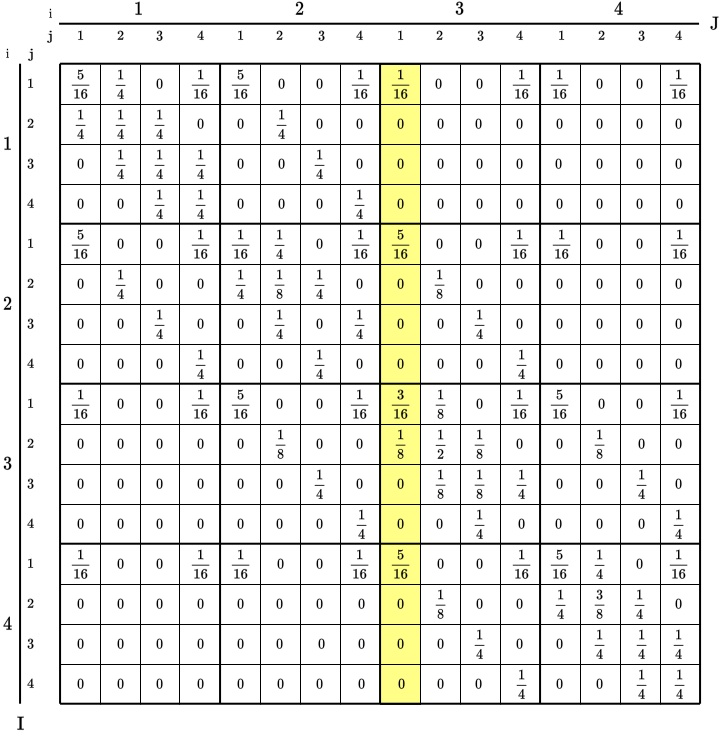}}
	\caption{Example of the definition of the transition matrix $\mathsf{M}$. The probabilities that a particle at position $J=(3,1)$ moves to a position $I=(i,j)$ are for $I=(1,1):\frac{1}{4} \cdot \frac{1}{4}$, for $I=(2,1)$ and $I=(4,1):\frac{1}{4} + \frac{1}{4} \cdot \frac{1}{4}$, for $I=(3,2):\frac{1}{8}$ and $\frac{3}{16}$ for the remaining position $J = (3,1)$ (left). The corresponding complete transition matrix $\mathsf{M}$ with highlighted column $J=(3,1)$ (right)}
	\label{fig:Transition_Matrix}
\end{Figure}
\section{UResNet: U-shaped Residual Convolutional Neural Network}\label{UResNet}
U-shaped neural networks (U-Nets) have been introduced by \cite{RonFisBro:2015:ucn} for applications in medical image segmentation. However, these networks are increasingly employed for applications in different fields. For instance, pixel-wise regression within the proposed U-Net enables the modeling of magnetic stray fields. For a more efficient training a residual U-Net is presented as residual neural network propose skip connections or shortcuts across some layers. 
\subsubsection*{Deep Convolutional Neural Network Architecture}
To predict the components of the stray magnetic fields $\widetilde{\bH}$ per pixel, a modified version of the U-shaped residual neural network model based on \cite{TakNas:2020:add} is created. The detailed architecture is shown in Fig. \ref{fig:UResNet_Architecture}. The input of the model is the two-dimensional binary image data, in which a pixel with a value of 0 (1) corresponds to the magnetic solid (vacuum space). From the UResNet model, we obtain a magnetic scattering field as $128 \times 128$ output with the values of vector components per pixel. Generally, the desired pixel size of images in Deep Convolutional Neural Networks is $2^m$, in this case, $2^7$.
\begin{Figure}[htbp]
  \includegraphics[width=16.0cm]{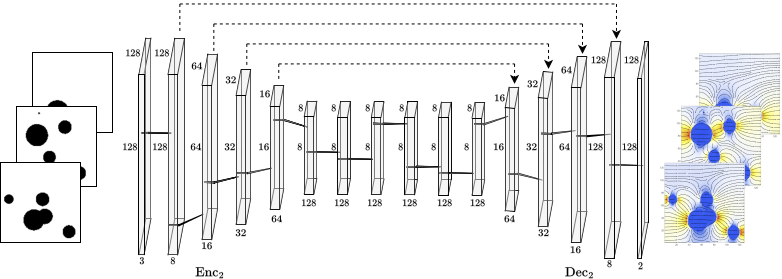}
  \caption{UResNet architecture set up}
  \label{fig:UResNet_Architecture}
\end{Figure} 

There are 5 downsampling blocks in the encoding part of the network, which are used to reduce the number of variables and compress the geometric information into a reduced dimension image.  
 Each block includes a pooling layer, two convolution layers, each followed by a batch normalization layer, a parametric rectified linear unit (PReLU) activation function.  
A detailed architecture of the second coding block is given in Tab. \ref{Structure_Enc_2}. Pooling, located in the first layer of each block, results in a reduction of the image size by half, thus reducing training time and computational cost. 
While there are several pooling methods, in this model the most prevalent max-pooling method was used, where only the maximum output value is retained for all subsequent layers and computations when reducing the dimensions. 
 Each convolutional layer has a number of output channels, also called filters, and a kernel of size $3 \times 3$. 
 The filters are used to extract important key features from the corresponding input data. 
  In convolution, the input data is multiplied by the matrices of the filter kernels, and the result is summed and combined into an output pixel. This process is repeated for the entire image, shifting the filter along the image with a stride of one.
  A commonly used regularization method in deep neural networks is the use of dropout layers at the end of each block to reduce the risk of overfitting. 
 However, batch normalization already regularizes the model and reduces the need for dropout, see \cite{SriHinKri:2014:das}.
From investigations with and without dropouts, we discovered that removing dropouts entirely results in higher validation accuracy for the network in this case.

\begin{Table}[htb!]
\centering
\caption{Detailed structure for downsampling encoding block $"\text{Enc}_2"$}
\label{Structure_Enc_2}
\begin{tabular}{lccc}
\toprule
\multicolumn{3}{l}{\textbf{Detailed structure for downsampling encoding block $"\textbf{Enc}_2"$}} \\
\midrule
Operation Layer    & Filters & Input Size & Output Size \\
\midrule
Pooling Layer & 8 & $128 \times 128 \times 8$  & $64 \times 64 \times 8$\\
Convolutional Layer	& 8 & $64 \times 64 \times 8$  & $64 \times 64 \times 8$ \\
Batch Normalization Layer	& 8 & $64 \times 64 \times 8$  & $64 \times 64 \times 8$ \\
Parametric Rectified Linear Unit (PReLU)	& 8 & $64 \times 46 \times 8$  & $64 \times 64 \times 8$\\
Convolutional Layer				& 16 & $64 \times 64 \times 8$  & $64 \times 32 \times 16$\\
Batch Normalization Layer		& 16 &  $64 \times 64 \times 16$  & $64 \times 64 \times 16$ \\
Parametric Rectified Linear Unit (PReLU)		& 16 & $64 \times 64 \times 16$  & $64 \times 64 \times 16$\\
\bottomrule
\end{tabular}
\end{Table}

Four residual blocks are located between the encoding and decoding blocks, which can improve gradient flow during training and simplify the network by effectively skipping connections. 
Unlike the traditional neural network structure, residual neural networks not only have all layers lined up but also skip connections, connecting blocks of equal size from downsampling and upsampling. 
The blocks have a mirrored layer structure due to a chosen U-shaped design. 
The five upsampling blocks for the decoding part of the network are designed to gradually scale the image to the pixel size of the original input. 
The first layer of each decoding block is concatenated with the corresponding coding block. 
Convolutional neural networks use multiple deconvolution layers in image resampling and formation, iteratively resampling a larger image from a set of lower resolution images. 
However, instead of a deconvolution layer, a resizing layer is used to resize the image by a factor of two using the nearest neighbor interpolation method to avoid checkerboard artifacts in the output images. 
Therefore, each decoding block specified in Tab. \ref{Structure_Dec_2} has a similar structure to the encoding block with two convolutional layers, each followed by a stack normalization layer and a parametric rectified linear unit activation function.
 Within the last layer, the input is mapped to the desired output size of $128 \times 128 \times 2$, which corresponds to the magnetic stray field $\widetilde{\bH}$.
\begin{Table}[htb!]
\centering
\caption{Detailed structure for upsampling decoding block "Dec$_2$"}
\label{Structure_Dec_2}
\begin{tabular}{lccc}
\toprule
\multicolumn{3}{l}{\textbf{Detailed structure for upsampling decoding block $"\textbf{Dec}_2"$}} \\
\midrule
Operation Layer    & Filters & Input Size & Output Size \\
\midrule
Convolutional Layer	& 16 & $ 64 \times 64 \times 16 $  & $64 \times 64 \times 16$ \\
Batch Normalization Layer	& 16 & $64 \times 64 \times 16$  & $64 \times 32 \times 16$ \\
Parametric Rectified Linear Unit (PReLU)	& 16 & $64 \times 64 \times 16$  & $64 \times 64 \times 16$\\
Convolutional Layer				& 8 & $64 \times 64 \times 16$  & $64 \times 32 \times 8$\\
Batch Normalization Layer		& 8 &  $64 \times 64 \times 8$  & $64 \times 64 \times 8$ \\
Parametric Rectified Linear Unit (PReLU)		& 8 & $64 \times 64 \times 8$  & $64 \times 64 \times 8$ \\
Resize Layer & 8 & $64 \times 64 \times 8$  & $128 \times 128 \times 8$\\
\bottomrule
\end{tabular}
\end{Table}

\section{Fourier Convolutional Neural Network}\label{FCNN}
The second approach is a Fourier Convolutional Neural Network (FCNN), first introduced by \cite{MatHenLeC:2014:fto,PraWilCoe:2017:Ffc}, which provides a simple alternative to traditional DCNNs by computing convolutions in the Fourier domain instead of the spatial domain. Convolution in the spatial domain corresponds to an element-wise product in the Fourier domain, which can significantly reduce the computational cost of training.
\subsubsection*{Fourier Convolutional Neural Network Architecture}
The Fourier Convolutional Neural Network consists of two fully connected layers (FC) in the first and last layers respectively, and four consecutive Fourier layers.
In this model the architecture setup of the Fourier Neural Operator (FNO), introduced by \cite{LiKovAzi:2020:fno}, was adapted for the application on the prediction of magnetic stray fields. The network architecture is illustrated in Fig. \ref{fig:FNO_Architecture}. 
\begin{Figure}[htbp]
  \includegraphics[width=16.0cm]{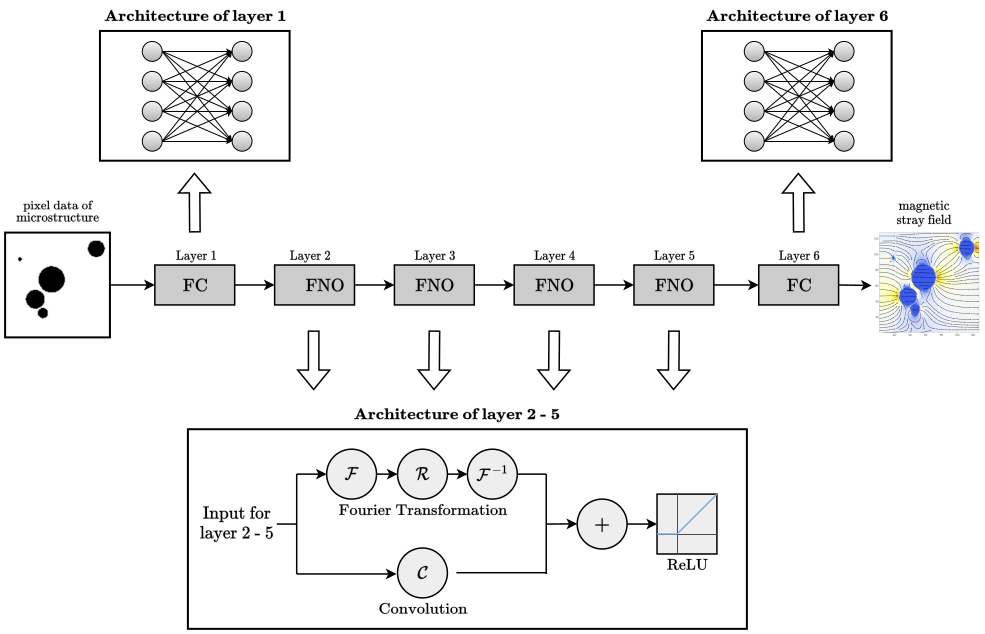}
  \caption{Fourier Convolutional Neural Network (FCNN) architecture}
  \label{fig:FNO_Architecture}
\end{Figure}

For the application of the Fourier Neural Operator numerical inputs are required. Hence, the binary structures of all images in the data set are encoded with prior data preprocessing, transforming the input into a 2-D matrix of pixel data instead of image input. Although The output of the magnetic stray fields is also a 2-D matrix, in this case, no further decoding is required as the output is already in the desired form.  
Based on the preprocessed pixel data of the binary image, the network starts by transforming the input into higher dimensions using a Fully Connected (FC) Layer with a number of 32 channels. Here, the channels represent the width of the FNO network, generally the number of features to be learned in each layer.
In the top path, the input is encoded by applying the discrete Fourier transform (DFT), defined as
\begin{equation}
    \mathcal{F}(u,v) = \sqrt{\frac{2}{N}} \sqrt{\frac{2}{M}}  \sum_{i=0}^{N-1} \sum_{j=0}^{M-1}{ \mathcal{F}^{-1}(i,j) \exp \left(-2\pi\iota \, \left(\frac{u\,i}{N}+\frac{v\,j}{M}\right)\right)}
\end{equation}
A linear transformation $\mathcal{R}$ is applied in the spectral domain to filter multiple high frequencies, since high frequencies often contain noise and imply abrupt changes in the image. 
In this context, the Modi parameter specifies the degree of signal filtering. In other words, the modes define the number of lower Fourier modes retained. Here, a number of 12 Fourier modes is selected.
Noting that the maximum allowed number of retained modes depends on the size of the computer grid, see \cite{GuaHsuChi:2021:fno}.
The hyperparameter setting, namely the chosen number of channels and modes, is one of the most effective methods to adjust the accuracy of the model.
While convolution is performed in frequency space in the top path, the bottom path is implemented solely in spatial domain. In the latter, the input is transformed by a one-dimensional spatial convolution $\mathcal{C}$ with a kernel of size 1 and a stride of 1.
The two paths are recombined after decoding the frequency domain values into the spatial dimension using the inverse Discrete Fourier Transform, which is obtained by
\begin{equation}
    \mathcal{F}^{-1}(i,j) = \sqrt{\frac{2}{N}} \sqrt{\frac{2}{M}}  \sum_{u=0}^{N-1} \sum_{v=0}^{M-1}{ \mathcal{F}(u,v) \exp \left(2\pi\iota \, \left(\frac{u\,i}{N}+\frac{v\,j}{M}\right)\right)}
\end{equation}
where $\iota$ with $\iota^2=-1$ is the imaginary unit. 
Once the two paths are added, a non-linear activation function, the Rectified Linear Unit activation, is applied. 
Finally, the received feature maps are resized to the desired dimensions using a fully connected layer.
Fourier Neural Operators can model the complex operators in partial differential equations, which often exhibit highly non-linear behavior and high frequencies, by the interaction of linear convolution, Fourier transformations, and the non-linear activation function, see \cite{GuaHsuChi:2021:fno}.
\section{Numerical Examples}\label{sec:Example}
For the evaluation of the representative capability and accuracy of the surrogate models presented in section \ref{UResNet} and \ref{FCNN}, the UResNet and the FCNN, three different aspects are considered: the analysis of the training process and loss, the performance evaluation on a new example and the comparison of the correlation. At first, a detailed look at the training process is given and different regression loss metrics are evaluated. Then, the predictive performance is illustrated on a new example from the unseen test set. At last, a comparison of the correlation of both models is given for the different subsets, indicating the best and worst results. The solution obtained by Brownian motion simulations is referred to as ground truth in the following.
\subsection{Training Process and Regression Loss Metrics}
To assess the quality of the predictive modeling that involved predicting a numeric value for the magnetic stray field, we compare some regression loss metrics. 
Three error measures are typically used for analyzing and reporting a model's performance: The Mean Squared Error (MSE), the Mean Absolute Error (MAE), and the Rooted Mean Squared Error (RMSE). The sum of the squared differences between the expected components of the magnetic stray field \mbox{\small $\widetilde{\bH}$} and the predicted magnetic stray field \mbox{\small $\widetilde{\bH^{\ast}}$} by the models is the MSE. In contrast, the MAE is the absolute deviation, and the square root of the MSE denotes the RMSE. In addition, we also consider relative deviations from prediction and ground truth in the error metrics given by the normalized MSE (NMSE) and normalized RMSE (NRMSE).
The aforementioned metrics for loss assessment are now considered and  later compared for each training pair $\tau$, which are defined as
\vspace{-1mm}
\begin{align*}
\centering
\textrm{Mean Squared Error:} &&
\textrm{$f^{\tau}$ := MSE}_{\tau} &= \frac{1}{n}\sum_{i=1}^n{\frac{1}{m}\sum_{j=1}^m{\left({\widetilde{H^{\ast}_{ij}} - \widetilde{H_{ij}}}\right)^2}} \\
\textrm{Normalized MSE:} &&
\textrm{NMSE}_{\tau} &= \textrm{MSE}_{\tau} \Big/ \,{\overline{\widetilde{H}}^{\,2}} \\
\textrm{Mean Absolute  Error:} &&
\textrm{MAE}_{\tau} &= \frac{1}{n}\sum_{i=1}^n{\frac{1}{m}\sum_{j=1}^m{|{\widetilde{H^{\ast}_{ij}} - \widetilde{H_{ij}}}|}}\\
\textrm{Rooted MSE:} &&
\textrm{RMSE}_{\tau} &= \sqrt{\frac{1}{n}\sum_{i=1}^n{\frac{1}{m}\sum_{j=1}^m{\left({\widetilde{H^{\ast}_{ij}} - \widetilde{H_{ij}}}\right)^2}}} \\
\textrm{Normalized RMSE:} &&
\textrm{NRMSE}_{\tau} &= \textrm{RMSE}_{\tau} \Big/ \,{\overline{\widetilde{H}}} 
\end{align*}

with the mean value \mbox{\small ${\overline{\widetilde{H}}} = \frac{1}{n} \sum_{i=1}^n{||\widetilde{\bH_i}||} \neq 0$ }. 
Hereby, \mbox{\small ${\widetilde{\bH^{\ast}_i}}$} is the predicted and \mbox{\small $\widetilde{\bH_i}$} is the desired vector of the magnetic stray field per pixel. 
Before comparing the different error measures, the training process and error development is examined.
%\textcolor{red}{
For the implemented network, the final training process runs on a GPU (GeForce RTX 3090) and the update formula ADAM, see \cite{KinBa:2017:aam},
%shown in Tab. \ref{Tab:ADAM},
and a batch size of 8 is used for optimization.
%}
From an entire dataset of 500 examples, we randomly selected 400 training examples and 50 validation examples for training. In this case, the training time took about 1 minute, while predicting the trained network takes only about 0.012 seconds. An illustration of the relative error evolution is shown in Fig. \ref{fig:Loss_plot} over 100 training rounds, showing the rapidly decreasing error within the first training rounds for each of the presented models. 
The graph also depicts that UResNet model requires more rounds of iteration for the error to converge than the FCNN models. However, the final error only differs marginally and has similar order of magnitude.
From experiments with varying training rounds, it appeared that having a larger number of iteration rounds for the UResNet model would result in a further reduction in the relative error. For the FCNN models, an additional number of iterations would not cause the error to decrease further.

\begin{Figure}[htbp]
\centering
    \includegraphics[width=12.7cm]{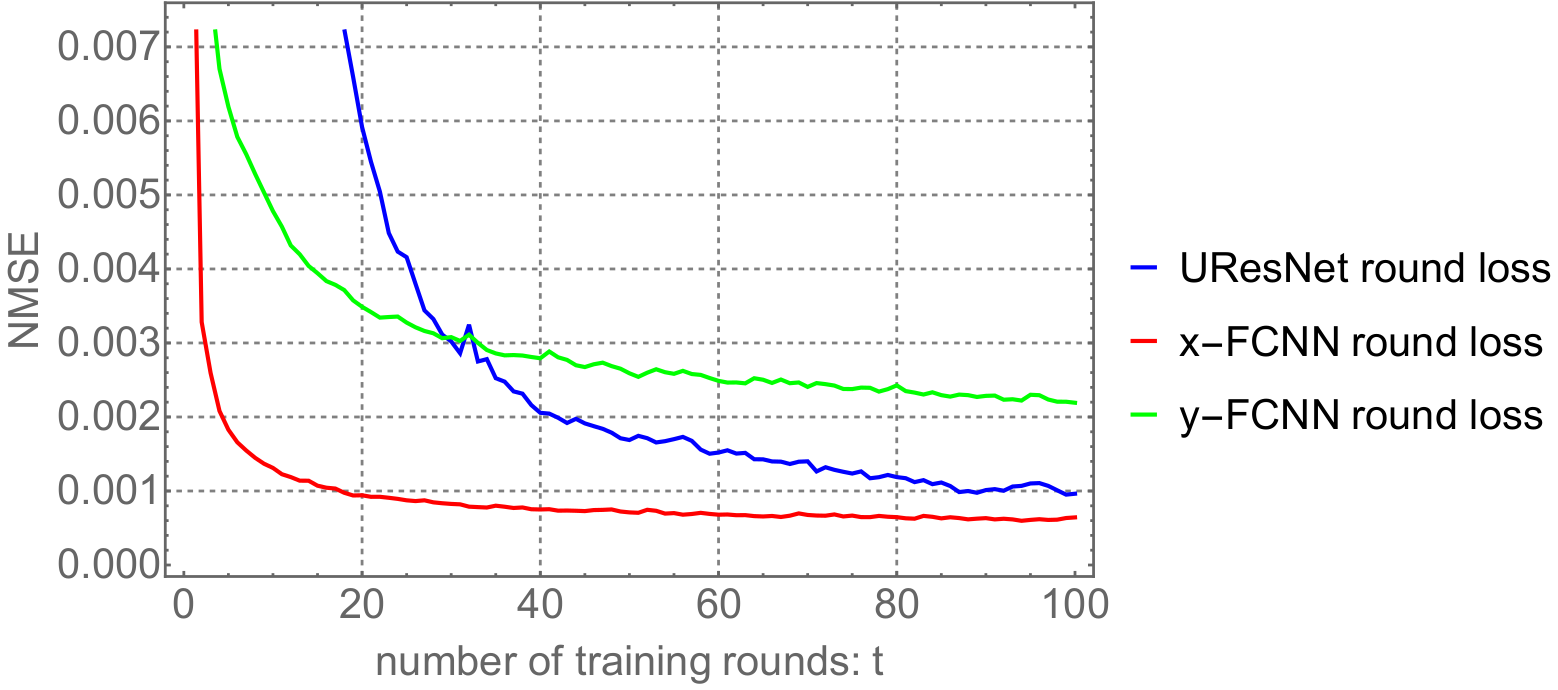}
 \caption{Loss plot of the arithmetic mean value of the NMSE during the training process for the whole training set, also known as round loss}
  \label{fig:Loss_plot}
 \end{Figure}
 
The two vector components of the UResNet model are trained concurrently, whereas the $x$ and $y$ components of the magnetic field are learned separately in the FCNN model. 
According to the computation, individual training of the components of the UResNet model resulted in no substantial increase in performance.
 Noticeably, there are more significant relative deviations for the $y$-component of the stray field than for the $x$-components. This deviation can be caused by the fact that the main direction of the total magnetic field is in the $x$-direction, resulting in the higher values and making it harder to learn the small values of the $y$-components. 
As shown in Tab. \ref{Tab:Regression_Loss_Metrics_Brown}, the trained models yield reasonable error measures for the training, validation, and test set. Each error measure was calculated for the prediction of the vector-valued magnetic stray field, i.e., \mbox{\small$\widetilde{\bH}$.}
\begin{Table}[htbp]
\centering
\caption{Values of different regression loss metrics after 100 training rounds using the surrogate UResNet and FCNN models}
\label{Tab:Regression_Loss_Metrics_Brown}
\begin{tabular}{cccccc}
\toprule
{\bf{UResNet}} & $\bf{MSE}$ & $\bf{NMSE}$ & $\bf{MAE}$ & $\bf{RMSE}$ & $\bf{NRMSE}$  \\
\midrule
round loss & 0.0009 & 0.0010 & 0.0198 & 0.0309 & 0.0317 \\ 
validation loss & 0.0017 & 0.0018 & 0.0249 & 0.0423 & 0.0435 \\
test loss & 0.0020 & 0.0021 & 0.0258 & 0.0449 & 0.0462\\[1mm]
\toprule
{\bf{FCNN}} & $\bf{MSE}$ & $\bf{NMSE}$ & $\bf{MAE}$ & $\bf{RMSE}$ & $\bf{NRMSE}$  \\
\midrule
round loss & 0.0004 & 0.0005 & 0.0123 & 0.0210 & 0.0217 \\ 
validation loss & 0.0013 & 0.0013 & 0.0196 & 0.0358 & 0.0360 \\
test loss & 0.0014 & 0.0015 & 0.0193 & 0.0377 & 0.0393 \\[1mm]
\bottomrule
\end{tabular}
\end{Table}

Since the validation set is only used to evaluate the model's performance on new data, it is not used for the training process. A significant gap between training, validation and test errors would indicate an overfitting condition in which the model solely learns the training data to the best extent possible. 
 
\subsection{Results}

The predictive accuracy of our surrogate models when applied to new unknown data is assessed by examining the performance of our trained models in predicting magnetic stray fields using the same sample for the UResNet in section \ref{sec:test_set_UResNet} and for the FCNN in section \ref{sec:test_set_FCNN} FCNN.
Here, a sample from the test set was selected since the examples from the training set are known through optimization. The validation set, although unknown to the models, was used to calculate the out-of-sample error during the training process. 
The results of our developed and trained UResNet and FCNN models with image dimensions of $128 \times 128$ pixels show similar behavior.
%\textcolor{red}{
%For completeness, Appendix \ref{sec:Examples} additionally display a numerical example from the validation and training sets analogously. !!!!Appendix fehlt!!!!}

\subsubsection{Example: Test Set (UResNet)}\label{sec:test_set_UResNet}

Taking an example from the test set, a qualitative comparison of the predicted stay fields from the trained UResNet model in Fig. \ref{fig:UResNet_example_test}b shows a high agreement with the reference data from the simulations based on Brownian motion in Fig. \ref{fig:UResNet_example_test}c.
 Given that this was a new unseen example to the neural network, this indicates a strong performance even for net unknown data. Within the representation of the two plots, only the x-component \mbox{\small $\widetilde{H_1}$} is depicted in the contour plots. 
 The color scaling was set uniformly for both illustrations, with the boundaries assigned to the minimum and maximum values from either the prediction \mbox{\small $\widetilde{H_1^{\ast}}$} or the ground truth \mbox{\small $\widetilde{H_1}$}.
  The entire field with its vector quantities of \mbox{\small $\widetilde{\bH}$} is plotted in the stream plots. 

\begin{Figure}[htb!]
\begin{minipage}[t]{0.325\textwidth}
\mbox{\small \quad (a) Input geometry}
\includegraphics[height=4.8cm]{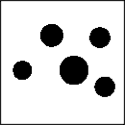}
\end{minipage}
\begin{minipage}[t]{0.325\textwidth}
\mbox{\small \, \, (b) UResNet prediction $\widetilde{\bH^{\ast}}$}
\includegraphics[height=5cm]{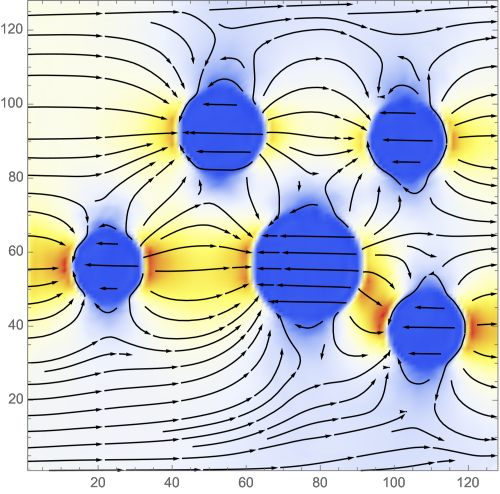}
\end{minipage}
\,
\begin{minipage}[t]{0.325\textwidth}
\mbox{\small $\qquad$ (c) Ground truth $\widetilde{\bH}$}
\includegraphics[height=5cm]{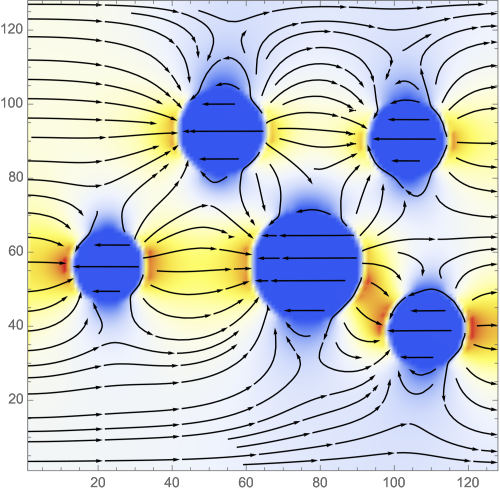}
\end{minipage}\\
\begin{minipage}[t]{0.9\textwidth}
\hspace{5.5cm} \includegraphics[height=1.5cm]{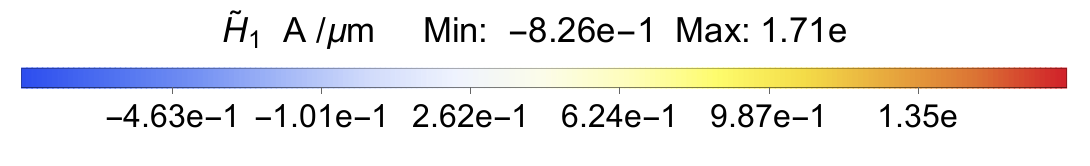}
\end{minipage}
\captionof{figure}{Magnetic stray fields for UResNet model with input dimensions $128\times128$. The color plots depict the x-component $\widetilde{H_1}$.}
\label{fig:UResNet_example_test}
\end{Figure}

 In addition to contour plots, the discrete values of the magnetic field $H_1$ and the magnetic stray field $\widetilde{H_1}$ along the horizontal intersection line A--A of the mid pixels are plotted in Fig. \ref{fig:A-A_plots_ex_test}. 
 The curve of the discrete predicted values (left) shows a similar behavior as the curve of the ground truth (middle), although it has a higher degree of oscillation. 
 This also emphasizes the good representative ability of the neural network. The distribution of the absolute error of the discrete values along the intersection line is displayed on the right using a histogram. The majority of errors are in the range of $0$ to $0.06$ with only a few individual "break-outs" and a maximum error of $0.14$. 

\begin{Figure}[htb!]  
\includegraphics[height=3.2cm]{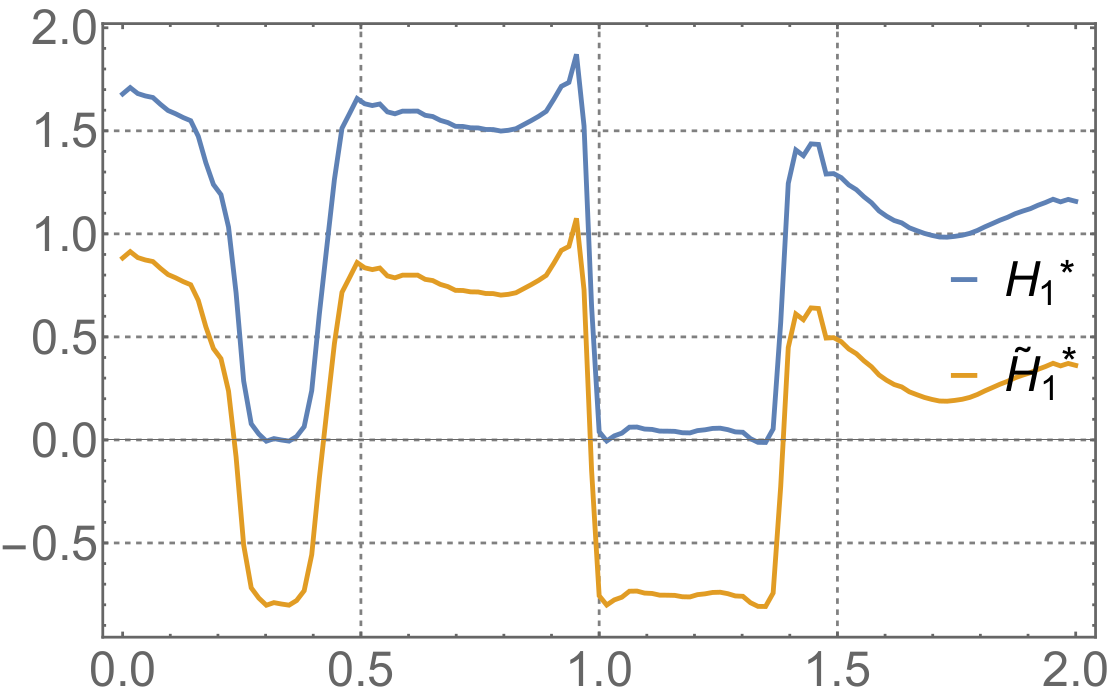}
\hspace{0.1cm}
\includegraphics[height=3.1cm]{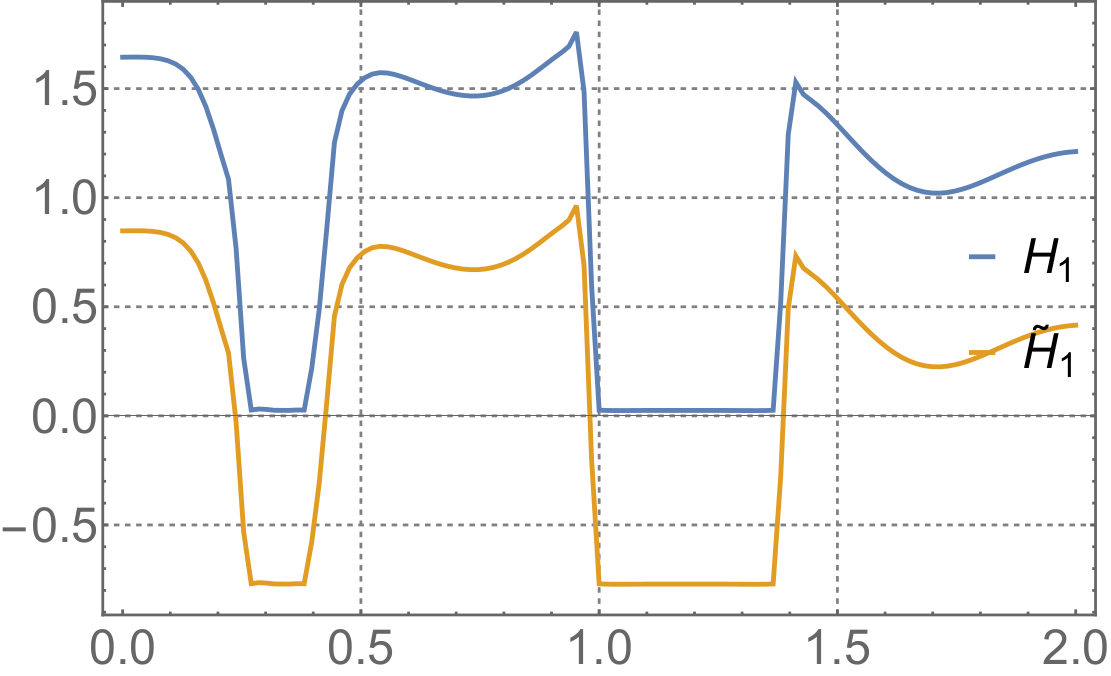}
\hspace{0.1cm}
\includegraphics[height=3.35cm]{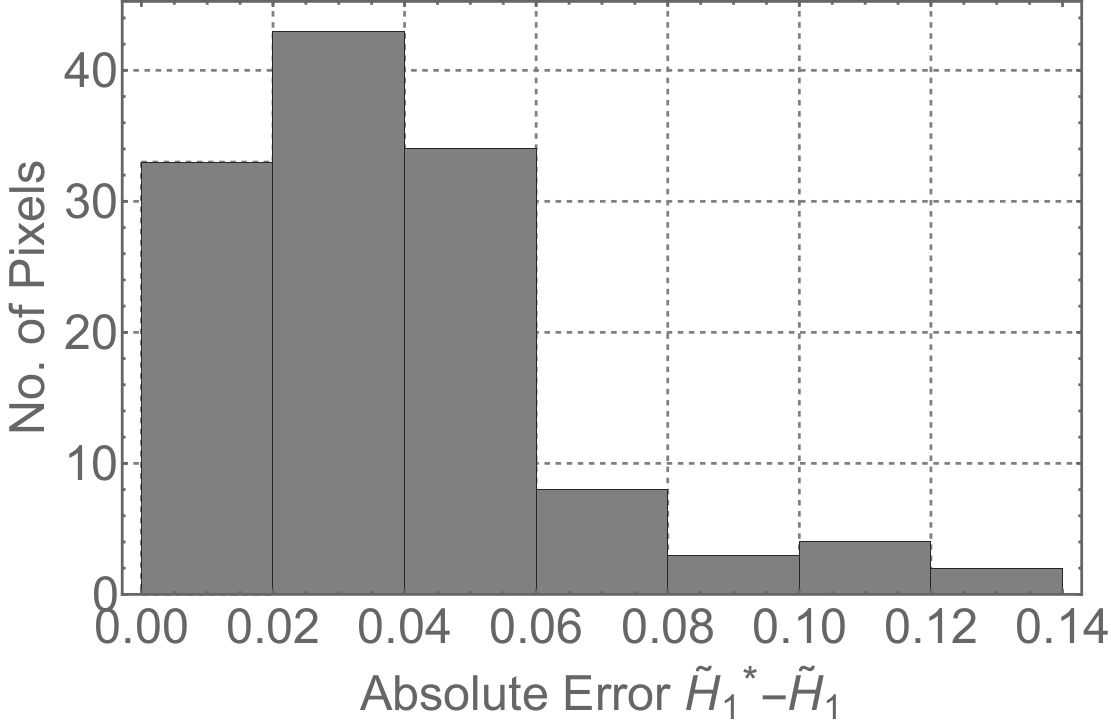}
\captionof{figure}{UResNet model: Discrete values of the total magnetic field $H_1$ and the fluctuation field $\widetilde{H_1}$ along the intersection line A--A are plotted for the prediction (left) and ground truth (middle) values. Distrubition of the absolute Error \mbox{\small $\widetilde{H_1^{\ast}} - \widetilde{H_1}$} per pixel (right).}
\label{fig:A-A_plots_ex_test}
\end{Figure}

The absolute error $\widetilde{\bH^{\ast}} - \widetilde{\bH}$ is shown in Fig. \ref{fig:Ex_test_error} where both plots represent the same numerical error with a different color scale. 
In the left plot, the automatic output of the contour plot of the absolute error is presented. However, the full data range is not used for color scaling. 
A so-called automatic plot range clipping is employed to get a visual idea of which pixels have the highest error, meaning the red regions do not correspond to the maximum error value and the blue regions do not necessarily correspond to the minimum error value. Values above $0.09$ are clipped to red, and values below $0.0085$ are clipped to blue. 
The absolute error over the entire data range is illustrated on the right. 
The bottom limit of the color scale is the minimum error value, and the upper limit is the maximal value of the absolute error for all $(128 \times 128)$ pixels of this example. 
Since the error is rather small in general and no red areas of high error values are visible, both plots are depicted for this example. 
 
\begin{Figure}[htbp]
\centering
    \includegraphics[height=5cm]{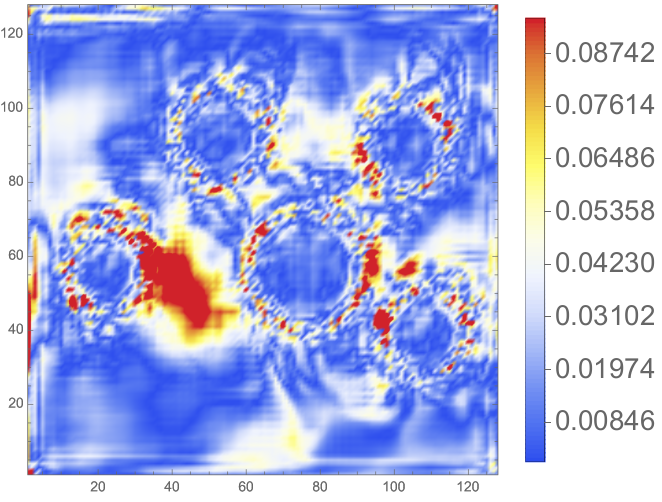}
    \hspace{0.5cm}
    \includegraphics[height=5cm]{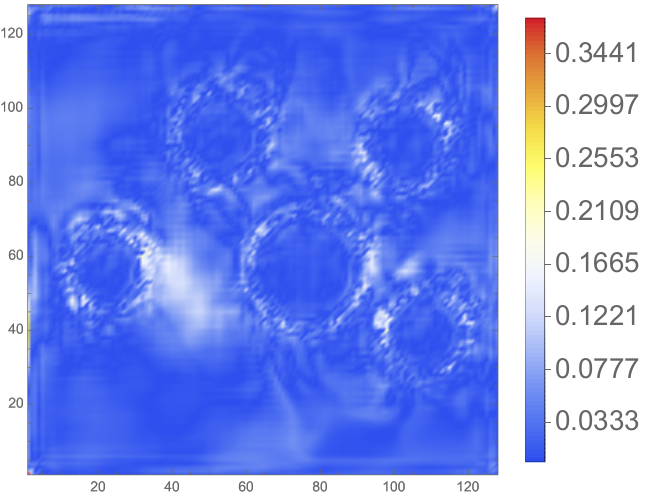}
 \caption{UResNet model: Illustration of the absolute error for a sample from the training set for image dimensions $128 \times 128$ with plot range clipping (left) and without plot range clipping (right).}
  \label{fig:Ex_test_error}
\end{Figure}

\subsubsection{Example: Test Set (FCNN)}\label{sec:test_set_FCNN}
The prediction for the unknown sample from the test set is also illustrated for the FCNN model. The side-by-side comparison of input geometry, FCNN prediction, and ground truth of Brownian motion is displayed in Fig. \ref{fig:FCNN_example_test}(a)--(c). Here, the second surrogate model also shows strong predictive performance on data withheld from the model entirely.

\begin{Figure}[htb!]
\begin{minipage}[t]{0.325\textwidth}
\mbox{\small \quad (a) Input geometry}
\includegraphics[height=4.8cm]{figures/Ex_test_img.png}
\end{minipage}
\begin{minipage}[t]{0.325\textwidth}
\mbox{\small \, \, (b) FCNN prediction $\widetilde{\bH^{\ast}}$}
\includegraphics[height=5cm]{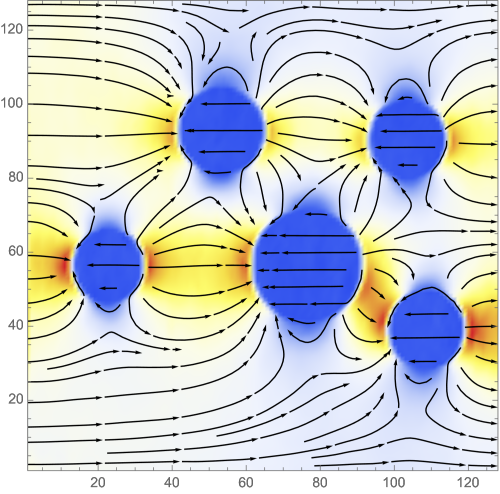}
\end{minipage}
\,
\begin{minipage}[t]{0.325\textwidth}
\mbox{\small $\qquad$ (c) Ground truth $\widetilde{\bH}$}
\includegraphics[height=5cm]{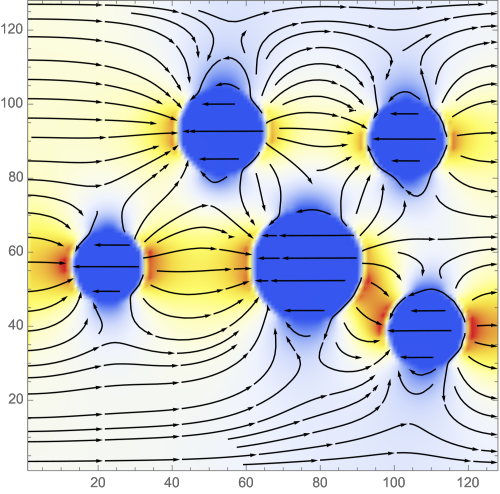}
\end{minipage}\\
\begin{minipage}[t]{0.9\textwidth}
\hspace{5.5cm} \includegraphics[height=1.5cm]{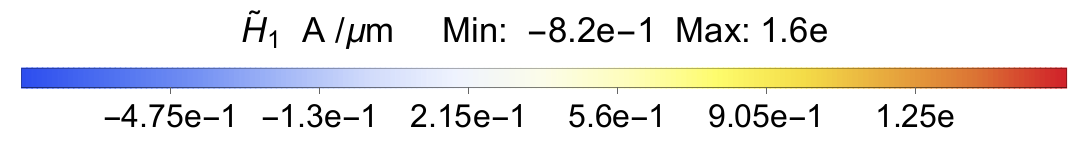}
\end{minipage}
\captionof{figure}{Magnetic stray fields for FCNN model with input dimensions $128\times128$. The color plots depict the x-component $\widetilde{H_1}$.}
\label{fig:FCNN_example_test}
\end{Figure}

The curves of the discrete values along the intersection line A--A in Fig. \ref{fig:A-A_plots_ex_test_FCNN} show a higher degree of oscillation for the predicted magnetic stray fields, similar to the UResNet model previously seen in Fig. \ref{fig:A-A_plots_ex_test}. 

\begin{Figure}[htb!]  
\includegraphics[height=3.1cm]{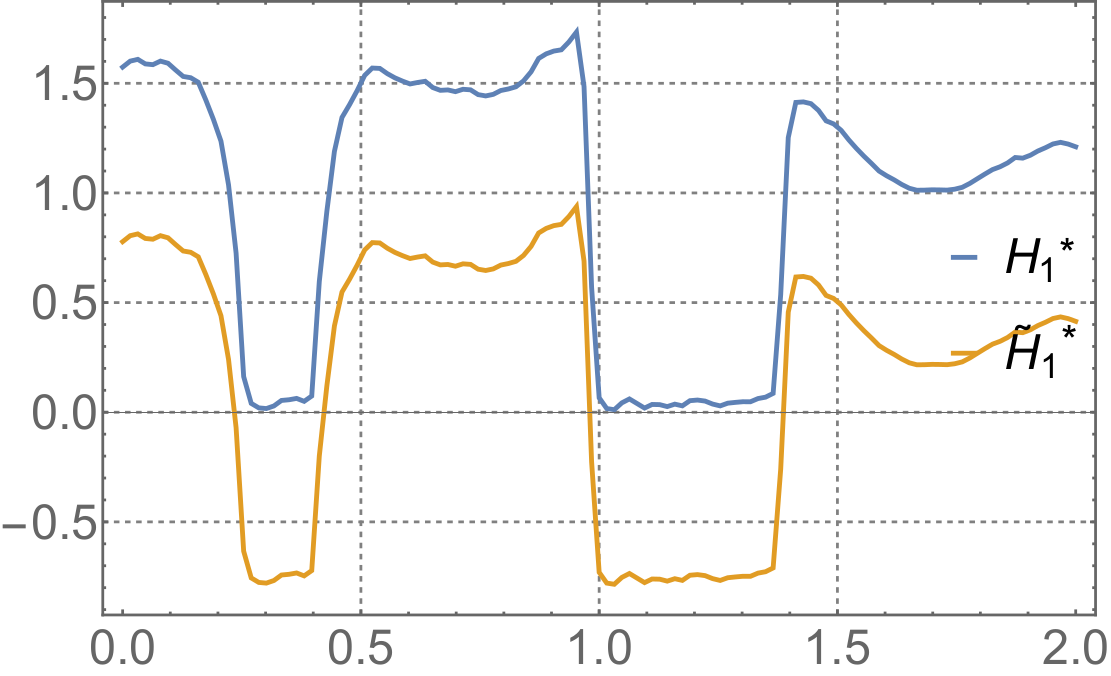}
\hspace{0.1cm}
\includegraphics[height=3.1cm]{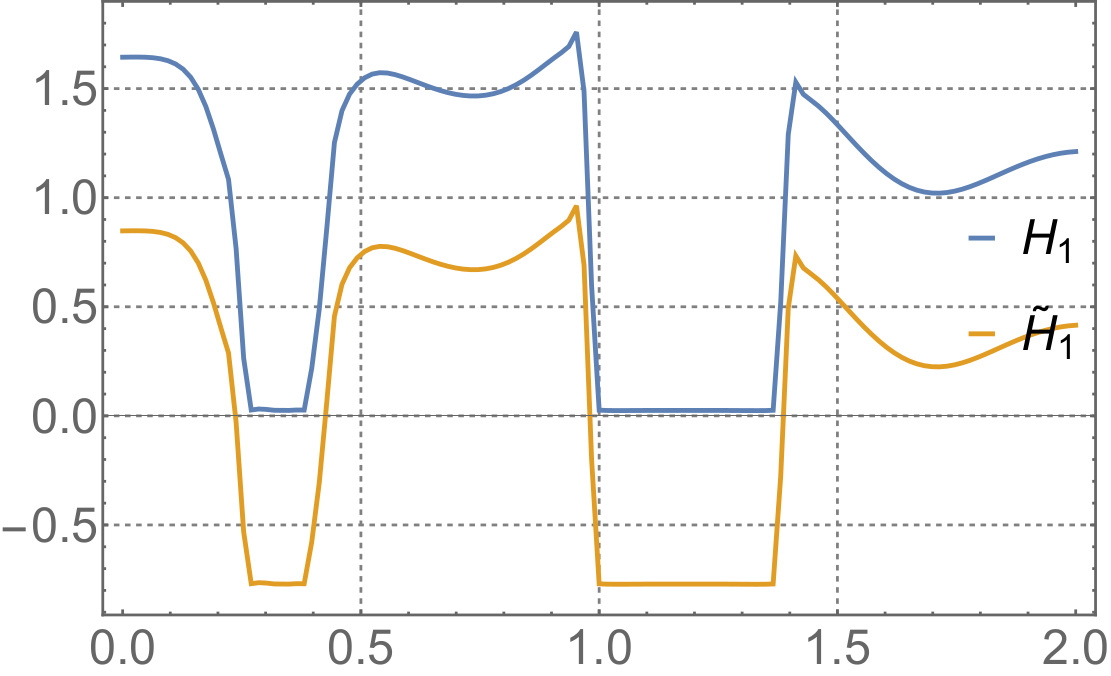}
\hspace{0.1cm}
\includegraphics[height=3.35cm]{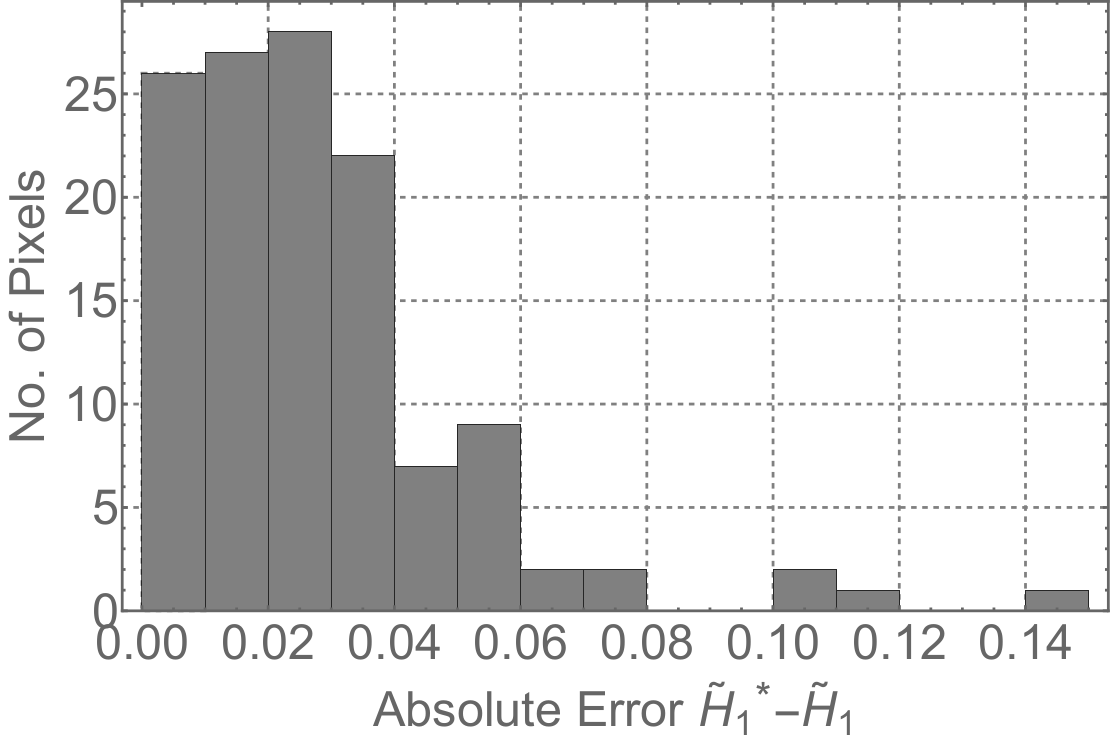}
\captionof{figure}{FCNN model: Discrete values of the total magnetic field $H_1$ and the fluctuation field $\widetilde{H_1}$ along the intersection line A--A are plotted for the prediction (left) and ground truth (middle) values. Distrubition of the absolute Error \mbox{\small $\widetilde{H_1^{\ast}} - \widetilde{H_1}$} per pixel (right).}
\label{fig:A-A_plots_ex_test_FCNN}
\end{Figure}

From the predictive comparison of the example from the test set, it may be challenging to evaluate which of the two surrogate models provides the more accurate prediction. However, looking at the absolute errors, Fig. \ref{fig:Ex_test_error_FCNN} shows that the maximum absolute error for the FCNN model prediction is indeed similar but slightly smaller. The maximum deviation for the predicted pixel value \mbox{\small $\widetilde{H_1^{\ast}}$} and the desired pixel value \mbox{\small $\widetilde{H_1}$} of ground truth is about 0.35 for the UResNet while having a value of about 0.32 for the FCNN.

\begin{Figure}[htbp]
\centering
    \includegraphics[height=5cm]{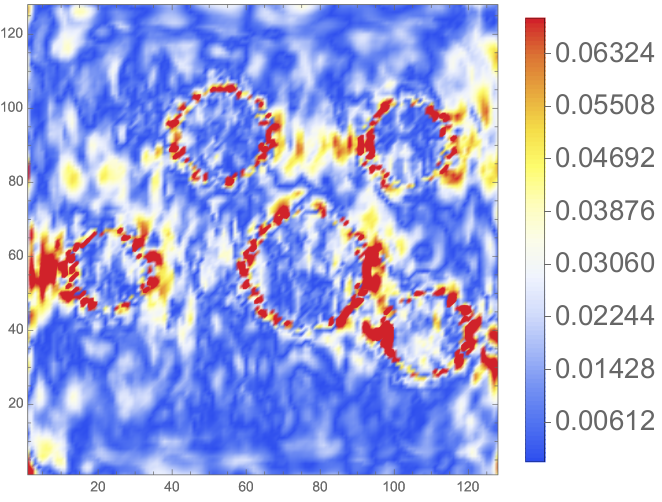}
    \hspace{0.5cm}
    \includegraphics[height=5cm]{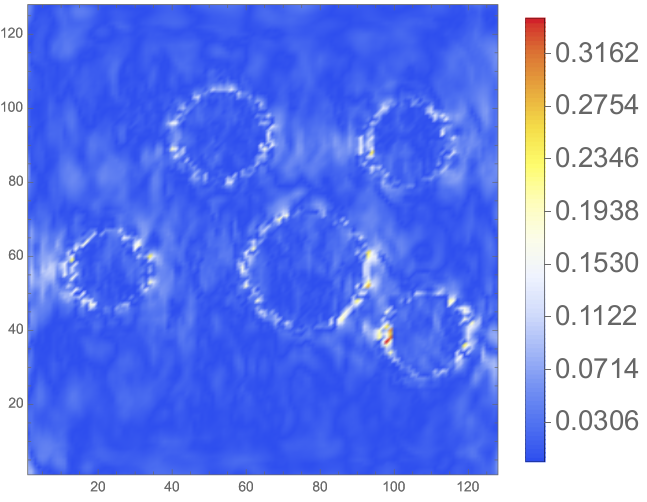}
 \caption{FCNN model: Illustration of the absolute error for a sample from the training set for image dimensions $128 \times 128$ with plot range clipping (left) and without plot range clipping (right).}
  \label{fig:Ex_test_error_FCNN}
\end{Figure}

The analysis of the results shown in Fig. \ref{fig:FCNN_example_test} -- \ref{fig:Ex_test_error_FCNN} were carried with the same plot settings as for the UResNet model for comparability reasons.

\subsection{Pearson Correlation Coefficient}
Correlation coefficients are one of the most commonly used metrics of analysis for the relationship between two pairs of data points, compare \cite{SchAxe:2011:cc}.
Here, the Pearson correlation coefficient (PCC) is used to calculate the correlation between the magnetic stray fields and their corresponding ground truths.
The larger the sum of the squares differences among the data, the lower the correlation coefficient.
Thus, the closer the predicted pixel values are to the desired solution, the higher the correlation coefficient.
The Pearson correlation coefficient is obtained from 
\begin{equation}
 \text{PCC} = \dfrac{\text{Cov}\left(\widetilde{\bH^{\ast}},\widetilde{\bH}\right)}{\sigma_{\widetilde{\bH^{\ast}}} \cdot \sigma_{\widetilde{\bH}}}
\end{equation}
with \mbox{\small $\sigma_{\widetilde{\bH^{\ast}}}$} and \mbox{\small $\sigma_{\widetilde{\bH}}$} being the standard deviations and Cov\mbox{\small$(\widetilde{\bH^{\ast}},\widetilde{\bH})$} as the covariance of the predicted values \mbox{\small$\widetilde{\bH^{\ast}}$} and the ground truth values \mbox{\small$\widetilde{\bH}$}. PCC scores of -1 represent entire disagreement while +1 represent total agreement. The correlation coefficient is 0 for entirely random guesses. 
A separate analysis of the correlation results from the different subsets is also given in Tab. \ref{Tab:Correlation_Metrics}, as well as the best and worst values.
The correlation results confirm the fact, already noticed during the error comparison that there are larger deviations in the predictions of the $y$-components. For each subset, the correlation coefficients of the $y$-conmponents of the stray field are slightly lower than for the $x$-components.

\begin{Table}[htbp]
\centering
\caption{Values of Pearson correlation coefficients for the different subsets of the trained surrogate models, indicating the worst and best correlations.}
\label{Tab:Correlation_Metrics}
\begin{tabular}{ccccccc}
\toprule
{\bf{UResNet}} & $\bf{PCC_{data}}$ & $\bf{PCC_{train}}$ & $\bf{PCC_{valid}}$ & $\bf{PCC_{test}}$ & $\bf{PCC_{worst}}$ & $\bf{PCC_{best}}$ \\[1mm]
\midrule 
$\widetilde{H_1}$ & $0.9965$ & $0.9971$ & $0.9947$ & $0.9932$ & $0.9838$ & $0.9986$ \\[1mm]
$\widetilde{H_2}$ & $0.9917$ & $0.9930$ & $0.9879$ & $0.9863$ & $0.9703$ & $0.9962$ \\[1mm]
$\widetilde{\bH}$ & $0.9984$ & $0.9987$ & $0.9976$ & $0.9971$ & $0.9813$ & $0.9999$ \\[1mm] 
\toprule
{\bf{FCNN}} & $\bf{PCC_{data}}$ & $\bf{PCC_{train}}$ & $\bf{PCC_{valid}}$ & $\bf{PCC_{test}}$ & $\bf{PCC_{worst}}$ & $\bf{PCC_{best}}$ \\[1mm]
\midrule
$\widetilde{H_1}$ & $0.9981$ & $0.9985$ & $0.9962$ & $0.9952$ & $0.9846$ & $0.9992$ \\[1mm] 
$\widetilde{H_2}$ & $0.9955$ & $0.9968$ & $0.9901$ & $0.9885$ & $0.9762$ & $0.9982$ \\[1mm] 
$\widetilde{\bH}$ & $0.9991$ & $0.9937$ & $0.9982$ & $0.9979$ & $0.9910$ & $0.9999$ \\[1mm]  
\toprule
\end{tabular}
\end{Table}

The square of the Pearson correlation coefficient is commonly referred to as the coefficient of determination, R$^2$, which is the proportion of the variation in one variable that can be explained by the variation in the other variables, see \cite{Kir:2008:pcc}. 
%\textcolor{red}{ The evolution of the coefficient of determination over 100 training rounds can be reported within the training progress measurements after each training round, as shown in Appendix \ref{sec:Training_Progress_Measurements} denoted in the column "res. R$^2$". !!!! Der Anhang fehlt!!!}

\begin{Figure}[htb!]
\begin{minipage}[t]{0.47\textwidth}
\centering
\includegraphics[width=\textwidth]{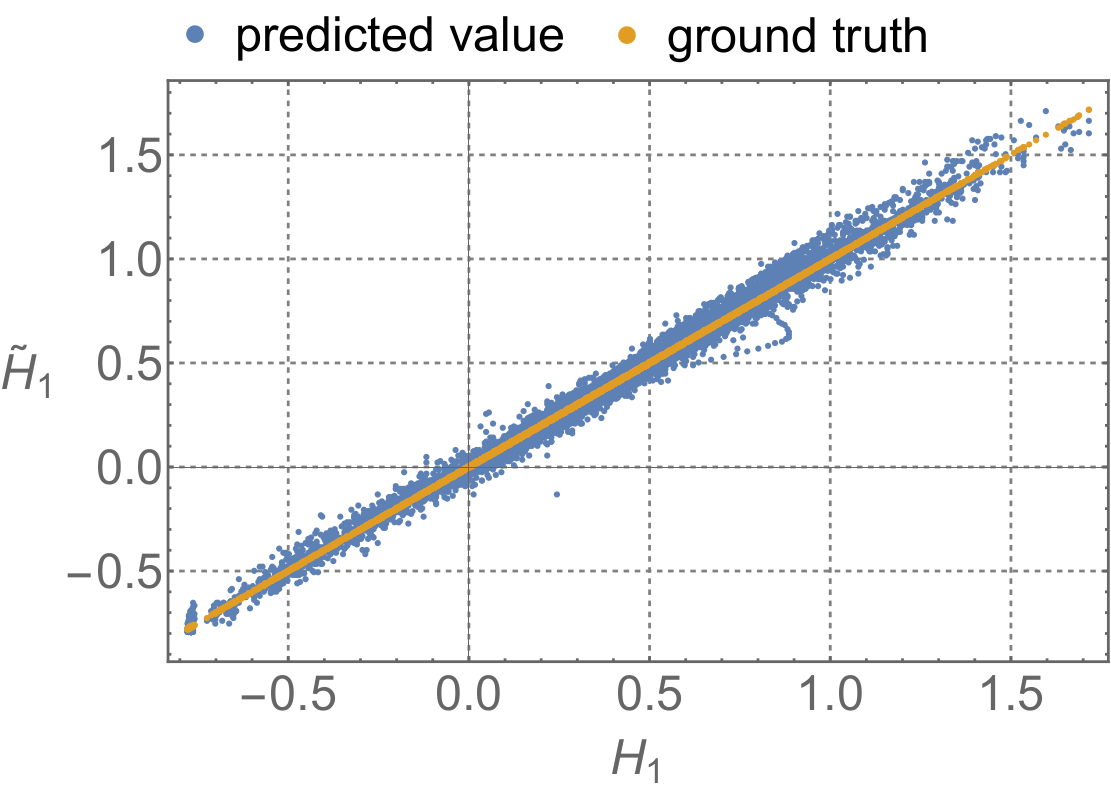} 
\end{minipage}
\begin{minipage}[t]{0.475\textwidth}
\hspace{0.5cm}
\includegraphics[width=\textwidth]{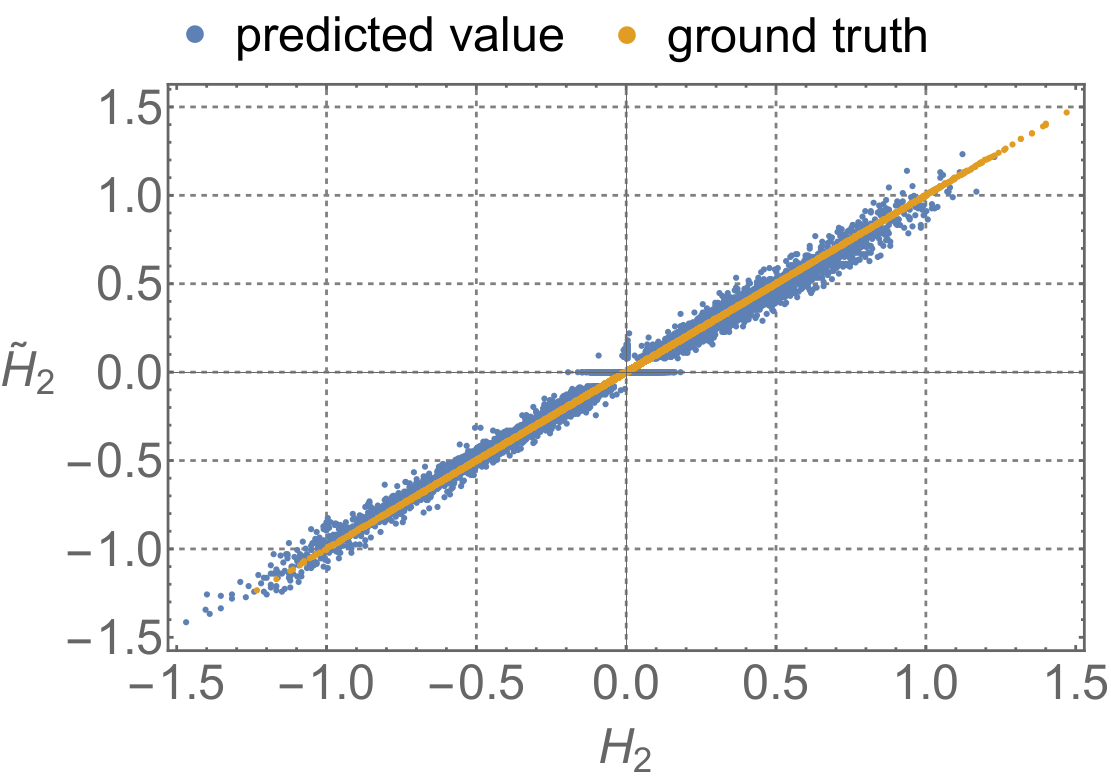} 
\end{minipage}
\captionof{figure}{UResNet: Cross plot of the predicted components $\widetilde{H^{\ast}}_1$ (left) and $\widetilde{H^{\ast}}_2$ (right) for each pixel vs. the ground truth components $\widetilde{H}_1$ and $\widetilde{H}_2$ of the system depicted in Fig. \ref{fig:UResNet_example_test}.}
\label{fig:x_y_Correlation_UResNet}
\end{Figure}

A correlation between two parameters can be examined visually with the help of a scatter plot of the two components.
For the sample in Fig. \ref{fig:UResNet_example_test}, the scatterplots of our predicted stray field for test images vs. ground truth in Fig. \ref{fig:x_y_Correlation_UResNet} indicate that the surrogate model has a high PCC of 0.9973 for $\widetilde{H}_1$ (left), a PCC of 0.9907 for $\widetilde{H}_2$ (right).
The yellow points display an ideal linear relationship, which corresponds to a correlation of the value $1$, while the blue points represent the predicted values vs. ground truth.

\begin{Figure}[htb!]
\begin{minipage}[t]{0.47\textwidth}
\centering
\includegraphics[width=\textwidth]{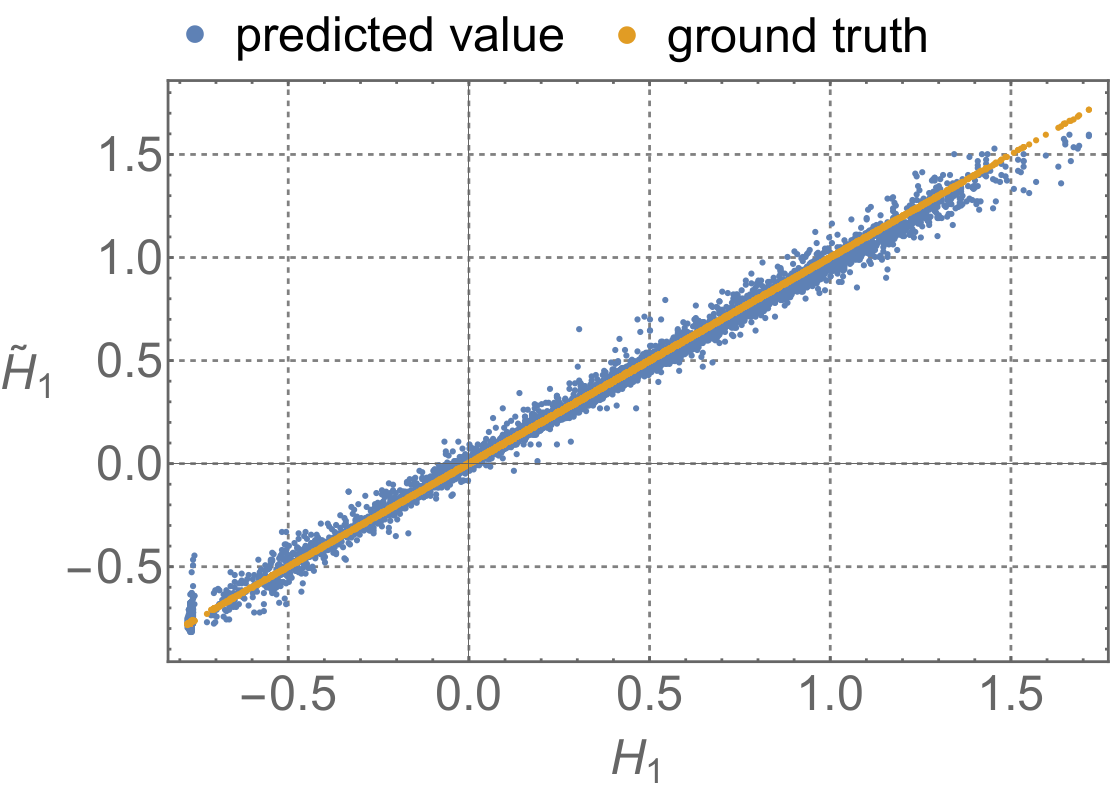} 
\end{minipage}
\begin{minipage}[t]{0.475\textwidth}
\hspace{0.5cm}
\includegraphics[width=\textwidth]{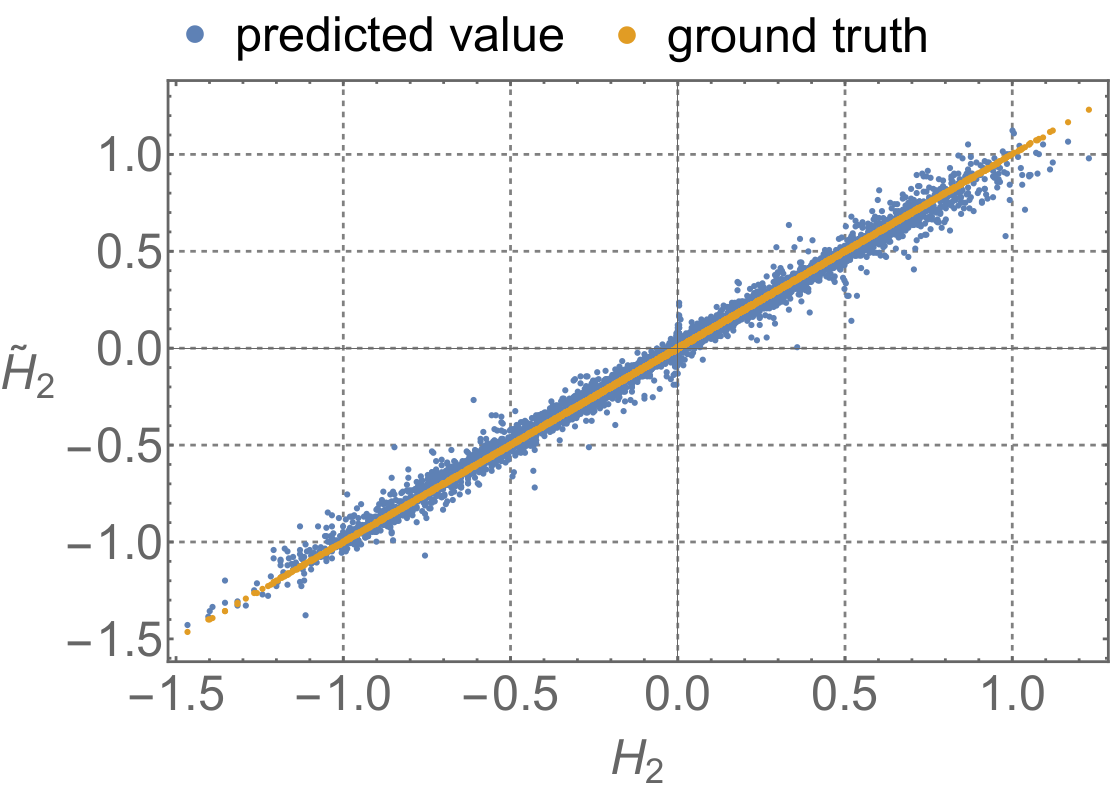} 
\end{minipage}
\captionof{figure}{FCNN: Cross plot of the predicted components $\hat{\widetilde{H}}_1$ (left) and $\hat{\widetilde{H}}_2$ (right) for each pixel vs. the ground truth components $\widetilde{H}_1$ and $\widetilde{H}_2$ of the system depicted in Fig. \ref{fig:FCNN_example_test}.}
\label{fig:x_y_Correlation_FCNN}
\end{Figure}

\vspace{0.3cm}

Fig. \ref{fig:x_y_Correlation_FCNN} shows the scatter plot of the FCNN model prediction for the sample in Fig. \ref{fig:FCNN_example_test} against ground truth. The correlation coefficients for the depicted plots yield even higher values of 0.9984 for $\widetilde{H}_1$ (left) and 0.9951 for $\widetilde{H}_2$ (right) than the UResNet prediction.

\section{Conclusion}
The presented work proposed two data-driven deep learning methods for calculating demagnetization and magnetic stray fields.
For the purpose of data generation, we use a low-cost stochastic model based on Brownian motion, which has the advantage of allowing the generation of a large number of data sets in a comparatively short time compared to traditional numerical methods such as the finite element method.
However, the surrogate models enable an even more efficient and faster prediction of magnetic stray fields with high accuracy than the simulations of stochastic Brownian motion and FEM.
Evaluations of the predictions on independent geometries were performed using a numerical example, showing a high correlation of the order of about $0.998$ and a relative error of $0.1-0.2$ \% (NMSE).
Here, the prediction performance of the trained FCNN model is slightly stronger than the trained UResNet model.
Since the presented surrogate models are applied in different scientific fields, they can also be used to approximate or predict other scientific scenarios that are not necessarily related to magnetic materials. 

\subsection*{Acknowledgement}
We gratefully acknowledge the financial support of the German Research Foundation (DFG) in the framework of the CRC/TRR 270, for 
project ZINF \glqq Management of data obtained in experimental and in silico investigations\grqq, 
project A07 \glqq Scale-bridging of magneto-mechanical mesostructures of additive manufactured and severe plastically deformed materials\grqq 
as well as for 
project A05 \glqq Designing 4f-3d permanent magnets by tailoring crystal fields\grqq, 
all with the project number 405553726.

% === list of references
% ------- layout-datei --------------
%\bibliographystyle{sn-aps}
\bibliographystyle{plainnat}
%\bibliographystyle{plaindin}
%\bibliographystyle{plaindin_shortname2}
%\bibliographystyle{elsarticle-num-names}
% ------- bib-datei --------------
%\bibliography{magneto}
\bibliography{ML_data_driven_04_2023_04_08.bib}

\end{document}